\documentclass{era-l}

\issueinfo{00}% volume number
  {}%         % issue number
  {}%         % month
  {1995}%     % year

\copyrightinfo{1995}%            % copyright year
  {American Mathematical Society}% copyright holder

\usepackage{pstricks}
\usepackage{amsmath}
\usepackage{latexsym}
\usepackage{mathrsfs}

\input{boxedeps}
\SetEPSFDirectory{} 
\HideDisplacementBoxes

\SetRokickiEPSFSpecial  %% dvips by Tom Rokicki

\theoremstyle{definition}

\theoremstyle{remark}

\numberwithin{equation}{section}

\DeclareMathAlphabet{\ams}{U}{msb}{m}{n}

\newpsobject{showgrid}{psgrid}{subgriddiv=1,griddots=5,gridlabels=6pt}

\def\R{\ams{R}}
\def\Z{\ams{Z}}

\def\CC{\mathscr C}
\def\LL{\mathscr L}

\def\aa{\alpha}

\def\o{\text{O}}
\def\reflec{\text{Reflec}}
\def\po{\text{PO}}

\def\vol{\text{vol}}
\def\covol{\text{covol}}
\def\isom{\text{Isom}\,}

\begin{document}

\title{The smallest hyperbolic 6-manifolds}

%    Information for first author
\author{Brent Everitt}
%    Address of record for the research reported here
\address{Department of Mathematics, University of York, York
YO10 5DD, England}
\email{bje1@york.ac.uk}
%    \thanks will become a 1st page footnote.
\thanks{The first author is grateful to the Mathematics Department, Vanderbilt University
for its hospitality during a stay when the results of this paper were obtained.}

%    Information for second author
\author{John Ratcliffe}
\address{Department of Mathematics, Vanderbilt University, Nashville, TN 37240, USA}
\email{ratclifj@math.vanderbilt.edu}

%    Information for third author
\author{Steven Tschantz}
\address{Department of Mathematics, Vanderbilt University, Nashville, TN 37240, USA}
\email{tschantz@math.vanderbilt.edu}

%    General info

\begin{abstract}
By gluing together copies of an all-right angled Coxeter polytope 
a number of open hyperbolic $6$-manifolds
with Euler characteristic $-1$ are constructed. They are the first known
examples of hyperbolic $6$-manifolds having the smallest possible volume. 
\end{abstract}

\maketitle

\section{Introduction}

The last few decades has seen a surge of activity in the study of
finite volume hyperbolic manifolds--that is, 
complete Riemannian $n$-manifolds of constant sectional curvature $-1$. 
Not surprisingly for geometrical objects, volume has been, and continues to be, the most 
important invariant for understanding their sociology.
The possible volumes in a fixed dimension forms a well-ordered subset of $\R$,
indeed a discrete subset except in $3$-dimensions (where the orientable manifolds at
least have ordinal 
type $\omega^\omega$). Thus it is a natural problem with a long history to construct
examples of manifolds with minimal volume in a given dimension.

In $2$-dimensions the solution is classical, with the minimum volume
in the compact orientable case achieved by a genus $2$ surface, and in the non-compact
orientable case by a once-punctured torus or thrice-punctured sphere
(the identities of the manifolds are of course also known in the non-orientable case). 
In $3$-dimensions the compact
orientable case remains an open problem with the Matveev-Fomenko-Weeks 
manifold \cite{Matveev-Fomenko88,Weeks85} obtained
via $(5,-2)$-Dehn surgery on the the sister of
the figure-eight
knot complement conjecturally the smallest. 
Amongst the non-compact orientable $3$-manifolds the figure-eight knot complement realizes
the minimum volume \cite{Meyerhoff01} and the Gieseking manifold (obtained by 
identifying the sides of a regular hyperbolic tetrahedron as in \cite{Everitt02b,Prok98})
does so for the non-orientable ones \cite{Adams87}.
One could also add ``arithmetic'' to our list of adjectives and so have eight
optimization problems to play with (so that %for example, 
the Matveev-Fomenko-Weeks manifold is known
to be the minimum volume orientable, arithmetic compact $3$-manifold, see \cite{Chinburg01}).

When $n\geq 4$ the picture is murkier, although in even dimensions we have
recourse to the Gauss-Bonnet Theorem, so that in particular the minimal volume
a $2m$-dimensional hyperbolic manifold could possibly have is when the 
Euler characteristic $\chi$ satisfies $|\chi|=1$. 
The first examples of non-compact $4$-manifolds with $\chi=1$ were constructed
in \cite{Ratcliffe00} (see also \cite{Everitt02}). The compact case remains a 
difficult unsolved problem, although if we restrict to arithmetic manifolds,
then it is known [2,8] that a minimal volume arithmetic compact orientable 4-manifold $M$
has $\chi \leq 16$ and $M$ is isometric to the orbit space of a torsion-free subgroup 
of the hyperbolic Coxeter group $[5,3,3,3]$. 
The smallest compact hyperbolic $4$-manifold currently known to exist 
has $\chi=8$ and is constructed in \cite{Conder04}.
Manifolds of very small volume have been constructed in $5$-dimensions 
\cite{Everitt02,Ratcliffe04},
but the smallest volume $6$-dimensional example hitherto known has $\chi=-16$ 
\cite{Everitt02}.

In this paper we announce the discovery of a number of non-compact non-orientable 
hyperbolic $6$-manifolds with Euler characteristic $\chi=-1$. The method of construction
is classical in that the manifolds are obtained by identifying the sides of a 
$6$-dimensional hyperbolic Coxeter polytope.

\section{Coxeter polytopes}

Let $C$ be a convex (not necessarily bounded) polytope of finite volume in a simply
connected space $X^n$ of constant curvature. Call $C$ a Coxeter polytope if the 
dihedral angle subtended by two intersecting $(n-1)$-dimensional sides 
is $\pi/m$ for some integer $m\geq 2$. When $X^n=S^n$ or the Euclidean space $E^n$,
such polyhedra have been completely classified \cite{Coxeter34}, 
but in the hyperbolic space 
$H^n$, a complete classification remains a difficult problem 
(see for example \cite{Vinberg85} and the references there).

If $\Gamma$ is the group generated by reflections in the $(n-1)$-dimensional sides
of $C$, then $\Gamma$ is a discrete cofinite subgroup of the Lie group $\isom X^n$,
and every discrete cofinite reflection group in $\isom X^n$ arises in this way 
from some Coxeter polytope, which is uniquely defined up to isometry. 
The Coxeter symbol for $C$ (or $\Gamma$) has nodes
indexed by the $(n-1)$-dimensional sides, and an edge labeled $m$ joining the nodes
corresponding to sides that intersect with angle $\pi/m$ (label the edge joining
the nodes of non-intersecting sides by $\infty$). In practice the labels $2$ and $3$
occur often, so that edges so labeled are respectively removed or left unlabeled.

Let $\Lambda$ be a $(n+1)$-dimensional Lorentzian lattice, 
that is, an $(n+1)$-dimensional 
free $\Z$-module equipped with a $\Z$-valued bilinear form of signature $(n,1)$. 
For each $n$, there is a unique such $\Lambda$, denoted
$I_{n,1}$, that is odd and self-dual (see \cite[Theorem V.6]{Serre73},
or \cite{Milnor73, Neumaier83}).
By \cite{Borel62}, the
group $\o_{n,1}\Z$ of automorphisms of $I_{n,1}$ acts 
discretely, cofinitely by isometries on the hyperbolic space $H^n$ obtained by 
projectivising the negative norm vectors
in the Minkowski space-time $I_{n,1}\otimes\R$ (to get a faithful action one normally
passes to the centerless version $\po_{n,1}\Z$). 

Vinberg and Kaplinskaja showed \cite{Vinberg78,Vinberg72} that
the subgroup $\reflec_n$ of $\po_{n,1}\Z$ generated by reflections in positive norm vectors
has finite index if 
and only if $n\leq 19$, thus yielding a family of cofinite reflection groups and 
corresponding finite volume Coxeter polytopes
in the hyperbolic spaces $H^n$ for $2\leq n\leq 19$. 
Indeed, Conway and Sloane have shown \cite[Chapter 28]{Conway93} or 
\cite{Conway82},
that for $n\leq 19$ the quotient 
of $\po_{n,1}\Z$ by $\reflec_n$ is a subgroup of the automorphism group of the Leech lattice.
Borcherds \cite{Borcherds87} showed that the (non self-dual) even sublattice of $I_{21,1}$
also acts cofinitely, yielding the highest dimensional example known of a Coxeter group
acting cofinitely on hyperbolic space.

When $4\leq n\leq 9$ the group $\Gamma=\reflec_n$ has Coxeter symbol,
$$
\begin{pspicture}(0,0)(6,2)
%\showgrid
\rput(3,1){\BoxedEPSF{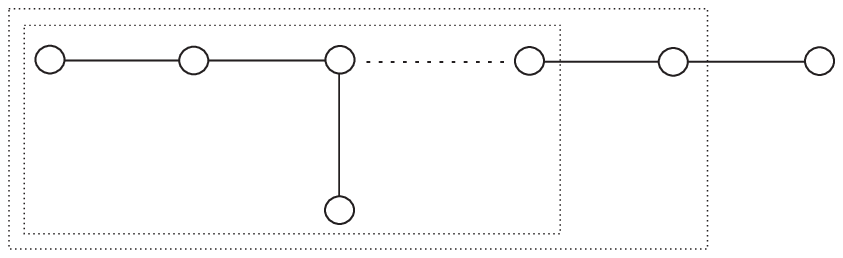 scaled 800}}
\rput(5.65,1.8){$4$}
\rput(5,.25){$\Gamma_v$}
\rput(3.8,.35){$\Gamma_e$}
\rput(6.25,1.2){$F_1$}
\rput(5.05,1.2){$F_2$}
\end{pspicture}
$$
with $n+1$ nodes and $C$ a non-compact, finite volume
$n$-simplex $\Delta^n$ (when $n>9$, the polytope $C$ has a more complicated 
structure).

Let $v$ be the vertex of $\Delta^n$ opposite the side $F_1$ marked on the symbol,
and let $\Gamma_v$ be the stabilizer in $\Gamma$ of this vertex. This stabilizer
is also a reflection group with symbol as shown, and is finite for 
$4\leq n\leq 8$ (being the Weyl group of type $A_4,D_5,E_6,E_7$ and $E_8$ respectively)
and infinite for $n=9$ (when it is the affine Weyl group of type $\widetilde{E}_8$). Let
$$
P_n=\bigcup_{\gamma\in\Gamma_v} \gamma(\Delta^n),
$$
a convex polytope obtained by gluing $|\Gamma_v|$ copies of the simplex $\Delta^n$ together. 
Thus, $P_n$ has finite volume precisely when $4\leq n\leq 8$, although it is non-compact, 
with a mixture of finite vertices in $H^n$ and cusped ones on $\partial H^n$.
In any case, $P_n$ is an all right-angled Coxeter
polytope: its sides meet with dihedral angle $\pi/2$ or are disjoint.
This follows immediately from the observation that the sides of $P_n$ arise from the
$\Gamma_v$-images of the side of $\Delta^n$ opposite $v$, and this side intersects the 
other sides of $\Delta^n$ in dihedral angles $\pi/2$ or $\pi/4$. 
Vinberg has conjectured that $n=8$
is the highest dimension in which finite volume all-right angled polytopes exist in 
hyperbolic space.

The volume of the polytope $P_n$ is given by
$$\vol(P_n)=|\Gamma_v|\vol(\Delta^n)=|\Gamma_v|[\po_{n,1}\Z:\Gamma]\covol(\po_{n,1}\Z),$$
where $\covol(\po_{n,1}\Z)$ is the volume of a fundamental region for the action
of $\po_{n,1}\Z$ on $H^n$ (and for $4\leq n\leq 9$ the index $[\po_{n,1}\Z:\Gamma]=1$). 
When $n$ is even, we have by \cite{Siegel36} and 
\cite{Ratcliffe97},
$$
\covol(\po_{n,1}\Z)=\frac{(2^{\frac{n}{2}}\pm 1)\pi^{\frac{n}{2}}}{n!}
\prod_{k=1}^{\frac{n}{2}} |B_{2k}|,
$$
with $B_{2k}$ the $2k$-th Bernoulli number and with the plus sign if 
$n\equiv 0,2\mod 8$ and the minus sign otherwise. 

Alternatively (when $n$ is even), we have
recourse to the Gauss-Bonnet Theorem, so that $\vol(P_n)=\kappa_n|\Gamma_v|
\chi(\Gamma)$, 
where $\chi(\Gamma)$ is the Euler characteristic of
the Coxeter group $\Gamma$ and $\kappa_n=2^n (n!)^{-1} (-\pi)^{n/2} (n/2)!$.
The Euler characteristic of Coxeter groups can be easily computed from their symbol
(see \cite{Chiswell92,Chiswell76} or \cite[Theorem 9]{Everitt02}). Indeed, when
$n=6$, 
$\chi(\Gamma)=-1/\LL$ where $\LL=2^{10}\,3^4\,5$ and so 
$\vol(P_6)=8\pi^3|E_6|/15\LL=\pi^3/15$.

The Coxeter symbol for $P_n$ has a nice description in terms of finite reflection groups.
If $v'$ is the vertex of $\Delta^n$ opposite the side $F_2$, let $\Gamma_e$ be the 
pointwise stabilizer of $\{v,v'\}$: the elements thus stabilize the edge $e$ of 
$\Delta^n$ joining $v$ and $v'$. 

Now consider the Cayley graph $\CC_v$ for 
$\Gamma_v$ with respect to the generating reflections in the sides of the symbol
for $\Gamma_v$. Thus, $\CC_v$ has vertices in one to one correspondence with the elements
of $\Gamma_v$ and for each generating reflection $s_{\aa}$, an undirected edge 
labeled $s_{\aa}$ connecting vertices $\gamma_1$ and
$\gamma_2$ if and only if $\gamma_2=\gamma_1s_{\aa}$ in $\Gamma_v$. In particular,
$\CC_v$ has $s_2$ labeled edges corresponding to the reflection in $F_2$. Removing these
$s_2$-edges
decomposes $\CC_v$ into components, each of which is a copy of the Cayley graph $\CC_e$ for
$\Gamma_e$, with respect to the generating reflections. 

Take as the nodes of the symbol for $P_n$ these connected components. If two components
have an $s_2$-labeled edge running between any two of their vertices in $\CC_v$, then leave 
the corresponding nodes
unconnected, otherwise, connect them by an edge labeled $\infty$. The resulting symbol
(respectively the polytope $P_n$) thus has $|\Gamma_v|/|\Gamma_e|$ nodes (resp. sides). 
The number of sides of $P_n$ for $n=4,5,6,7,8$ is $10,16,27,56$ and $240$ respectively. 

\section{Constructing the manifolds}
	
We now restrict our attention to the case $n=6$. 
We work in the hyperboloid model of hyperbolic 6-space
$$H^6=\{x\in \R^7: x_1^2+x_2^2+\cdots+x_6^2-x_7^2=-1\ \hbox{and}\ x_7>0\}$$
and represent the isometries of $H^6$ by Lorentzian $7\times 7$ matrices that 
preserve $H^6$. 
The right-angled polytope $P_6$ has 27 sides each congruent to $P_5$.  
We position $P_6$ in $H^6$ so that 6 of its sides are bounded 
by the 6 coordinate hyperplanes $x_i=0$ for $i=1,\ldots, 6$  
and these 6 sides intersect at the center $e_7$ of $H^6$. 
Let $K_6$ be the group of 64 diagonal Lorentzian $7\times 7$ matrices 
${\rm diag}(\pm 1,\ldots,\pm 1,1)$. 
The set $Q_6=K_6P_6$, which is the union of 64 copies of $P_6$, 
is a right-angled convex polytope with 252 sides. 
We construct hyperbolic 6-manifolds, with $\chi =-8$,
by gluing together the sides of $Q_6$ by a proper side-pairing 
with side-pairing maps of the form $rk$ with $k$ in $K_6$ and $r$ 
a reflection in a side $S$ of $Q_6$. 
The side-pairing map $rk$ pairs the side $S'=kS$ to $S$ 
(see \S 11.1 and \S 11.2 of \cite{Ratcliffe94} 
for a discussion of proper side-pairings). 
We call such a side-pairing of $Q_6$ simple. 
We searched for simple side-pairings of $Q_6$ that yield 
a hyperbolic 6-manifold $M$ with a freely acting $\Z/8$ symmetry group 
that permutes the 64 copies of $P_6$ making up $M$ in such a way 
that the resulting quotient manifold is obtained by gluing 
together 8 copies of $P_6$. 
Such a quotient manifold has $\chi = -8/8=-1$. 
This is easier said than done, since the search space 
of all possible side-pairings of $Q_6$ is very large. 
We succeeded in finding desired side-pairings of $Q_6$ by employing a strategy 
that greatly reduces the search space. 
The strategy is to extend a side-pairing in dimension 5 with the desired 
properties to a side-pairing in dimension 6 with the desired properties. 

Let $Q_5 = \{x\in Q_6: x_1 = 0\}$. 
Then $Q_5$ is a right-angled convex 5-dimensional polytope with 72 sides. 
Note that $Q_5$ is the union $K_5P_5$ of 32 copies of $P_5$ where 
$P_5= \{x\in P_6: x_1=0\}$ and  
$K_5$ is the group of 32 diagonal Lorentzian $7\times 7$ matrices 
${\rm diag}(1,\pm 1,\ldots,\pm 1,1)$.  
A simple side-pairing of $Q_6$ that yields a hyperbolic 6-manifold $M$ 
restricts to a simple side-pairing of $Q_5$ 
that yields a hyperbolic 5-manifold which is a totally geodesic hypersurface of $M$. 
All the orientable hyperbolic 5-manifolds that are obtained by gluing together 
the sides of $Q_5$ by a simple side-pairing are classified in \cite{Ratcliffe04}. 

We started with the hyperbolic 5-manifold $N$, numbered 27 in \cite{Ratcliffe04},  
obtained by gluing together the sides of $Q_5$ by the simple side-pairing 
with side-pairing code {\tt 2B7JB47JG81}. 
The manifold $N$ has a freely acting $\Z/8$ symmetry group 
that permutes the 32 copies of $P_5$ making up $N$ in such a way that 
the resulting quotient manifold is obtained by gluing together 4 copies of $P_5$. 
A generator of the $\Z/8$ symmetry group of $N$ is represented by the Lorentzian $6\times 6$ matrix
$$
\left(\begin{array}{cccccc}
     \phantom{-}1 & 0 &\phantom{-}0 & \phantom{-}1 & 0 &     -1 \\
     \phantom{-}0 & 0 &\phantom{-}0 & \phantom{-}0 & 1 & \phantom{-}0 \\
         -1 & 0 &          -1 & \phantom{-}0 & 0 & \phantom{-}1 \\
     \phantom{-}0 & 1 &\phantom{-}0 & \phantom{-}0 & 0 & \phantom{-}0 \\
     \phantom{-}0 & 0 &          -1 &     -1 & 0 & \phantom{-}1 \\
         -1 & 0 &          -1 &     -1 & 0 & \phantom{-}2
        \end{array} \right).
$$
The strategy is to search for simple side-pairings of $Q_6$ that yield 
a hyperbolic 6-manifold with a freely acting $\Z/8$ 
symmetry group with generator represented by the following Lorentzian $7\times 7$ matrix 
that extends the above Lorentzian $6\times 6$ matrix. 
$$
\left(\begin{array}{ccccccc}
     \phantom{-}1 &\phantom{-}0 & 0 &\phantom{-}0 & \phantom{-}0 & 0 &\phantom{-}0 \\
     \phantom{-}0 &\phantom{-}1 & 0 &\phantom{-}0 & \phantom{-}1 & 0 &     -1 \\
     \phantom{-}0 &\phantom{-}0 & 0 &\phantom{-}0 & \phantom{-}0 & 1 & \phantom{-}0 \\
          \phantom{-}0 &-1 & 0 &          -1 & \phantom{-}0 & 0 & \phantom{-}1 \\
     \phantom{-}0 &\phantom{-}0 & 1 &\phantom{-}0 & \phantom{-}0 & 0 & \phantom{-}0 \\
     \phantom{-}0 &\phantom{-}0 & 0 &          -1 &     -1 & 0 & \phantom{-}1 \\
         \phantom{-}0 &-1 & 0 &          -1 &     -1 & 0 & \phantom{-}2
        \end{array} \right)
$$      
For such a side-pairing the resulting quotient manifold can be obtained 
by gluing together 8 copies of $P_6$ by a proper side-pairing. 
By a computer search we found 14 proper side-pairings of 8 copies of
$P_6$ in this way, and hence we found 14 hyperbolic $6$-manifolds with $\chi = -1$. 
Each of these 14 manifolds is noncompact with 
volume $8\vol(P_6)=8\pi^3/15$ and five cusps. 
These 14 hyperbolic $6$-manifolds represent at least 7 different isometry types, 
since they represent 7 different homology types. 
Table 1 lists side-pairing codes for 7 simple side-pairings of $Q_6$ 
whose $\Z/8$ quotient manifold has homology groups isomorphic to  
$\Z^a\oplus(\Z/2)^b\oplus(\Z/4)^c\oplus(\Z/8)^d$ for nonnegative integers $a,b,c,d$
encoded by $abcd$ in the table. In particular, all 7 manifolds in Table 1 
have a finite first homology group.

All of our examples, with $\chi = -1$, can be realized as the orbit space $H^6/\Gamma$ 
of a torsion-free subgroup $\Gamma$ of $\po_{6,1}\Z$ of minimal index. 
These manifolds are the first examples of hyperbolic 6-manifolds 
having the smallest possible volume. 
All these manifolds are nonorientable. 
In the near future, we hope to construct orientable examples
of noncompact hyperbolic 6-manifolds having $\chi= -1$.

\begin{table}
\begin{center}
\begin{tabular}{lllllll}
$N$&$SP$&\ \ $H_1$&\ \ $H_2$&\ \ $H_3$&\ \ $H_4$&\ \ $H_5$\\
&&$\phantom{\mathbb Z}$0248&$\phantom{\mathbb Z}$0248&$\phantom{\mathbb Z}$0248&
$\phantom{\mathbb Z}$0248&$\phantom{\mathbb Z}$0248\\
1&{\tt GW8dNEEdN4ZJO1k2l1PIY}&
\phantom{${\mathbb Z}$}0401&
\phantom{${\mathbb Z}$}1910&
\phantom{${\mathbb Z}$}4821&
\phantom{${\mathbb Z}$}1500&
\phantom{${\mathbb Z}$}0000\\
2&{\tt HX9dNFEcM5aKU6f3f6UKa}&
\phantom{${\mathbb Z}$}0401&
\phantom{${\mathbb Z}$}1810&
\phantom{${\mathbb Z}$}8710&
\phantom{${\mathbb Z}$}5500&
\phantom{${\mathbb Z}$}0000\\
3&{\tt HX9dNFEcM5YIO1l3l1OIY}&
\phantom{${\mathbb Z}$}0401&
\phantom{${\mathbb Z}$}2900&
\phantom{${\mathbb Z}$}7810&
\phantom{${\mathbb Z}$}4400&
\phantom{${\mathbb Z}$}1000\\
4&{\tt HX9dNFEcM5YIO6l3l6OIY}&
\phantom{${\mathbb Z}$}0401&
\phantom{${\mathbb Z}$}2800&
\phantom{${\mathbb Z}$}7910&
\phantom{${\mathbb Z}$}4400&
\phantom{${\mathbb Z}$}1000\\
5&{\tt HX9dNFEcM5YIOxl3lyOIY}&
\phantom{${\mathbb Z}$}0211&
\phantom{${\mathbb Z}$}2800&
\phantom{${\mathbb Z}$}4821&
\phantom{${\mathbb Z}$}1400&
\phantom{${\mathbb Z}$}1000\\
6&{\tt HX9dNFEcM5YIOyl3lxOIY}&
\phantom{${\mathbb Z}$}0211&
\phantom{${\mathbb Z}$}2800&
\phantom{${\mathbb Z}$}4930&
\phantom{${\mathbb Z}$}1400&
\phantom{${\mathbb Z}$}1000\\
7&{\tt HX9dNFEcM5aKUxf3fyUKa}&
\phantom{${\mathbb Z}$}0301&
\phantom{${\mathbb Z}$}1900&
\phantom{${\mathbb Z}$}5630&
\phantom{${\mathbb Z}$}2500&
\phantom{${\mathbb Z}$}0000\\
\end{tabular}
\end{center}
\caption{Side-pairing codes and homology groups of the seven examples}
\end{table}

\bibliographystyle{amsplain}

\begin{thebibliography}{10}

\bibitem{Adams87}
C Adams.
\newblock {\em The non-compact hyperbolic $3$-manifold of minimum volume}.
\newblock Proc. AMS, 100 (1987), 601-606.

\bibitem{Belolipetsky02}
M Belolipetsky.
\newblock {\em On volumes of arithmetic quotients of $SO(1,n)$}.
\newblock Ann. Scuola Norm. Sup. Pisa Cl. Sci. (to appear).

\bibitem{Borcherds87}
R Borcherds.
\newblock {\em Automorphism groups of Lorentzian lattices}.
\newblock J. Algebra 111(1) (1987), 133-153.

\bibitem{Borel62}
A Borel and Harish-Chandra.
\newblock {\em Arithmetic subgroups of algebraic groups}.
\newblock Ann. of Math., 75(3) (1962), 485-535.

\bibitem{Chinburg01}
T Chinburg, E Friedman, K N Jones, and A W Reid. 
\newblock {\em The arithmetic hyperbolic 3-manifold of smallest volume}. 
\newblock Ann. Scuola Norm. Sup. Pisa Cl. Sci. (4) 30 (2001), 1-40. 

\bibitem{Chiswell92} 
I M Chiswell.
\newblock {\em The Euler characteristic of graph products and Coxeter groups}.
\newblock in Discrete Groups and Geometry, W J Harvey and Colin
Maclachlan (Editors), London Math. Soc. Lect. Notes, 173 (1992), 36-46.

\bibitem{Chiswell76} 
I M Chiswell.
\newblock {\em Euler characteristics of groups}.
\newblock Math. Z., 147 (1976), 1--11.

\bibitem{Conder04} 
M D E Conder and C Maclachlan.
\newblock {\em Small volume compact hyperbolic $4$-manifolds}.
\newblock to appear, Proc. Amer. Math. Soc.

\bibitem{Conway93}
J Conway and N J A Sloane.
\newblock {\em Sphere Packings, Lattices and Groups}.
\newblock Second Edition, Springer 1993.

\bibitem{Conway82} 
J Conway and N J A Sloane. 
\newblock {\em Leech roots and Vinberg groups}.
\newblock Proc. Roy. Soc. London Ser. A 384 (no. 1787) (1982), 233-258.

\bibitem{Coxeter34} 
H S M Coxeter. 
\newblock {\em Discrete groups generated by reflections}.
\newblock Ann. of Math, 35(2) (1934), 588-621.

\bibitem{Davis85} 
M W Davis.
\newblock {\em A hyperbolic 4-manifold.}
\newblock Proc. Amer. Math. Soc., 93 (1985), 325-328.

\bibitem{Everitt02}
B Everitt.
\newblock {\em Coxeter groups and hyperbolic manifolds.}
\newblock Math. Ann. 330(1), (2004), 127-150.

\bibitem{Everitt02b}
B Everitt.
\newblock {\em $3$-Manifolds from Platonic solids.}
\newblock Top. App. 138 (2004), 253-263.

\bibitem{Ratcliffe05}
B Everitt, J Ratcliffe and S Tschantz.
\newblock {\em Arithmetic hyperbolic $6$-manifolds of smallest volume}
\newblock (in preparation).

\bibitem{Matveev-Fomenko88}
V Matveev and A Fomenko.
\newblock {\em Constant energy surfaces of Hamilton systems, enumeration
of three-dimensional manifolds in increasing order of complexity, and
computation of volumes of closed hyperbolic manifolds.}
\newblock Russian Math. Surveys, 43 (1988), 3-24.

\bibitem{Meyerhoff01}
C Cao and G R Meyerhoff.
\newblock {\em The orientable cusped hyperbolic $3$-manifolds of minimum
volume}.
\newblock Invent. Math. 146(3) (2001), 451-478.

\bibitem{Milnor73}
J Milnor and D Husemoller.
\newblock {\em Symmetric Bilinear Forms}.
\newblock Springer, 1973.

\bibitem{Neumaier83}
A Neumaier and J J Seidel.
\newblock {\em Discrete Hyperbolic Geometry}.
\newblock Combinatorica 3 (1983), 219-237.

\bibitem{Prok98}
I Prok.
\newblock {\em Classification of dodecahedral space forms}.
\newblock Beitr\"{a}ge Algebra Geom. 39(2) (1998), 497-515.

\bibitem{Ratcliffe94}
J Ratcliffe.
\newblock {\em Foundations of hyperbolic manifolds}.
\newblock Graduate Texts in Mathematics 149, Springer 1994.

\bibitem{Ratcliffe97} 
J Ratcliffe and S Tschantz.
\newblock {\em Volumes of integral congruence hyperbolic manifolds}.
\newblock J. Reine Angew. Math. 488 (1997), 55-78.

\bibitem{Ratcliffe00} 
J Ratcliffe and S Tschantz.
\newblock {\em The volume spectrum of hyperbolic 4-manifolds.}
\newblock Experiment. Math. 9 (2000), no. 1, 101-125.

\bibitem{Ratcliffe04} 
J Ratcliffe and S Tschantz.
\newblock {\em Integral congruence two hyperbolic $5$-manifolds.}
\newblock Geom. Dedicata 107 (2004), 187-209.

\bibitem{Serre73}
J-P Serre.
\newblock {\em A Course in Arithmetic}.
\newblock Graduate Texts in Mathematics 7, Springer 1973.

\bibitem{Siegel36} 
C L Siegel.
\newblock {\em \"{U}ber die analytische Theorie der quadratischen Formen II}.
\newblock Ann. Math. 37 (1936), 230-263.

\bibitem{Vinberg85}
E B Vinberg.
\newblock {\em Hyperbolic Reflection Groups.}
\newblock Uspekhi Mat. Nauk 40:1 (1985), 29-66,
\newblock $=$ Russian Math. Surveys 40:1 (1985), 31-75.

\bibitem{Vinberg78}
E B Vinberg and I M Kaplinskaja.
\newblock {\em On the groups $\o_{18,1}(\Z)$ and $\o_{19,1}(\Z)$.}
\newblock SMD 19(1) (1978), 194-197.
\newblock $=$ DAN 238(6) (1978), 1273-1275.

\bibitem{Vinberg72}
E B Vinberg.
\newblock {\em On the groups of units of certain quadratic forms.}
\newblock Sb. 16(1) (1972), 17-35,
\newblock $=$ Mat. Sb. 87(1) (1972), 18-36.

\bibitem{Weeks85}
J Weeks. 
\newblock {\em Hyperbolic structures on $3$-manifolds}.
\newblock PhD thesis, Princeton University.

\end{thebibliography}

\end{document}